\documentclass{madrid15} 

\usepackage{pgfplots}\pgfplotsset{compat=1.4}
\usepackage[normalem]{ulem} 

\usepackage{graphicx}
\usepackage{caption}
\usepackage{subcaption}

\usepackage{amsmath}
\usepackage{amssymb}
\usepackage{amsfonts}
\usepackage{epsfig}
\usepackage{graphicx}
\usepackage{multicol}
\usepackage{color}
\usepackage{crop}
\usepackage{natbib}
\usepackage{float}

\usepackage{algorithmic}
\usepackage{algorithm}

\newcommand{\argmin}{\operatornamewithlimits{arg min}}

\newcommand{\bq}{\begin{eqnarray}}
\newcommand{\eq}{\end{eqnarray}}

\renewcommand{\eqref}[1]{equation (\ref{eq:#1})}

\newcommand{\qed}{\nobreak \ifvmode \relax \else
      \ifdim\lastskip<1.5em \hskip-\lastskip
      \hskip1.5em plus0em minus0.5em \fi \nobreak
      \vrule height0.75em width0.5em depth0.25em\fi}
\renewcommand{\Pi}{\mbox{\LARGE$\pi$}}

\begin{document}



	\begin{center}
		\large{\bf{General Optimization Framework for Robust and Regularized 3D Full Waveform Inversion}}\\
		\bigskip
		\normalsize{\emph{Stephen Becker, Lior Horesh, Aleksandr Aravkin, and Sergiy Zhuk}}\\
	\end{center}
	
	
\section{Introduction} 
Scarcity of hydrocarbon resources and high exploration risks motivate the development of high fidelity algorithms and computationally viable approaches to exploratory geophysics. 
Whereas early approaches considered least-squares minimization, recent developments have emphasized the importance of robust formulations, as well as formulations 
that allow disciplined encoding of prior information into the inverse problem formulation. 
The cost of a more flexible optimization framework is a greater computational complexity, as 
least-squares optimization can be performed using straightforward methods 
(e.g., steepest descent, Gauss-Newton, L-BFGS), whilst incorporation of robust (non-smooth) penalties requires custom changes that may be difficult to implement in the context of a general seismic inversion workflow. 
In this study, we propose a generic, flexible optimization framework capable of incorporating a broad range of noise models, 
forward models, regularizers, and reparametrization transforms. 
This framework covers seamlessly robust noise models (such as Huber and Student's $t$), as well as sparse regularizers, projected constraints, and Total Variation regularization. 
The proposed framework is also expandable --- we explain the adjustments that are required for any new formulation to be included. 
Lastly, we conclude with few numerical examples demonstrating the versatility of the formulation.


\section{Method and Theory}


\subsection{Classic Waveform Inversion} 
The canonical waveform inversion problem is to find the medium parameters for which the modeled data 
matches the recorded data in a least-squares sense \citep{tarantola84}. 
For linear wave equation  
\bq
H(m)u &=& q,
\eq
where $H(m)$ is the discretized wave equation of model parameters $m$, $u$ is the discretized wavefields and $q$ are the discretized source functions. 
The data are then given by sampling the wavefield at the receiver locations: $d=Su + \varepsilon$, where $\varepsilon$ stands for model misspecification and measurement noise. 
For example, for the simplest case of constant-density acoustics in the frequency domain, the wave equation is represented by the Helmholtz operator 
$H(m) = \omega^2 m + \nabla^2$ amended with appropriate boundary conditions, with squared slowness $m$.

The classic inverse problem can be cast as an explicit least-squares problem as follows:
\bq
\label{classicLS}
\min_{m} \;\phi(m) = \sum_{\omega}\|\underbrace{SH(m)^{-1}Q}_{h(m)} - D\|^2_F,
\eq
where $D = [d_1, d_2, \ldots d_N]$ 
corresponds to the recorded data vectors and 
$Q = [q_1,q_2,\ldots,q_N]$ are the corresponding source functions. 
$||\cdot||_F$ denotes the Frobenius norm which is defined as $||A||_F = \sqrt{\sum_{ij} a_{ij}^2}$. 

\subsection{Generalized Formulations}

To accommodate robust formulations, regularizers and constraints (e.g., Eq.~\ref{eq:rbox}), we consider 
\begin{equation}
\label{eq:gen}
\min_{y} \;\phi(y) = \sum_{\omega}\rho(h(Cy), D) + R(y),
\end{equation}
where $\rho$ is a general misfit penalty, and $R$ is a regularization function, and $C$ is a linear transformation to a space of 
interest (e.g. Fourier, wavelet, or curvelet). Thus, the classical inversion expression (Eq.~\ref{classicLS}) is a specific case where one set $R=0$, $\rho$ as $\rho_{LS}(\hat{D},D) = \|\hat{D}-D\|_F^2$, and $C=I$ (the identity).

We provide two notable examples of $R$.  
In the $1^{st}$ example, we encode \emph{constraints} using an indicator function $R$ that is \emph{infinite valued}.  
Let $B$ be a closed set of feasible values of the parameters $y$, and consider 
\begin{equation}
\label{eq:rbox}
R(y) = \begin{cases} 0 & \mbox{if} \quad y \in B\\
\infty & \mbox{if} \quad y \not\in B\end{cases}.
\end{equation}
The regularizer $R$ in this case is equivalent to a constraint $y\in B$, but we choose to express it as in Equation~\ref{eq:gen} 
because it offers a simple and uniform way to describe our algorithmic framework. 

For the second example, take $R$ to be $R_{lasso}(y)=\lambda \|y\|_1$. This is a sparsifying penalty on the transform space coefficients $y$; the larger 
the $\lambda$, the faster elements of $y$ are driven to $0$. In this case, $R$ is finite valued, but {\it not smooth}.  
We will also use $\ell_1$ \emph{constraints}, $R_{\ell_1-ball}$ defined by Eq.~\ref{eq:rbox} with $B(y)=\{ y : \|y\|_1 \le \tau \}$ for some parameter $\tau$. Projection onto this set can be performed in linear time.

The algorithmic framework we propose requires two assumptions:
\begin{enumerate}
\item The misfit penalty $\rho$ is differentiable. The differentiability assumption holds for least-squares, Huber~\citep{guitton2003robust}, hybrid~\citep{bube1997hybrid}, 
Student's $t$
~\citep{AravkinFHV:2012}, 
generalized self-tuning extensions~\citep{AravkinVanLeeuwen2012} amongst others. 
\item We require that the following problem can be solved quickly (e.g., linear time in length of $y$):
\[
\hat g := \argmin_g \frac{1}{2}\|g - y\|^2 + R(g).
\]
The vector $\hat g$ is known as the proximity, or {\it prox}, operator of $R$ applied to $y$, written as $\mbox{prox}_R(y)$.  
 
 \end{enumerate}

 The two examples we developed above
satisfy the $2^{nd}$ requirement for the $R$. For Eq.~\ref{eq:rbox},
 the prox operator is simply projection onto a set, 
which is often known in closed form. 
  For the non-smooth $R_{lasso}$ 
 the prox can also be implemented entry-wise, and is simply given by the soft thresholding operator 
 $(S_\lambda(y))_i = \text{sign}(y_i)\cdot \max(0,|y_i|-\lambda)$.
While a closed form formula may not be available for all cases, as long as $\mbox{prox}_R(y)$ can be efficiently computed, 
our algorithmic framework applies. 

\subsection{Algorithmic Framework}

To develop our algorithmic framework, we adopt and expand on the PQN algorithm of~\cite{SchmidtBergFriedlanderMurphy:2009}. 
This algorithm was developed for optimization of costly objective functions with simple constraints. Such settings are 
well suited to full-waveform inversion, since function evaluations require multiple expensive forward and adjoint PDE solves. PQN 
builds a quasi-Newton Hessian approximation for the objective from gradient information using a limited memory BFGS scheme~\citep{nocedal99}, 
and then solve an auxiliary subproblem
\[
\min_{y} Q(y) + R(y)
\]
where $Q$ is a limited memory-based convex quadratic, and $R$ represents an indicator function for simple constraints (such as Eq.~\ref{eq:rbox}).
The subproblem is solved by the spectral projected gradient method~\citep{birgin2000nonmonotone}, which relies upon assumption 2 we imposed in our framework, since prox is equivalent to projection when $R$ encodes constraints. 

In this study we expand this framework, allowing any function $R$ admitting an efficient proximity operator. In particular, we can then naturally 
incorporate sparsity promoting regularization, and total variation.

\begin{figure} 
\centering
\begin{tikzpicture}
  \begin{axis}[
    thick,
    width=.3\textwidth, height=.9cm, 
    xmin=-4,xmax=4,ymin=0,ymax=.5,
    no markers,
    samples=50,
    axis lines*=left, 
    axis lines*=middle, 
    scale only axis,
    xtick={-1,1},
    xticklabels={},
    ytick={0},
    ] 
\addplot[domain=-4:+4,densely dashed]{exp(-.5*x^2)/sqrt(2*pi)};
 \addplot[domain=-4:+4, dashdotted, red]{0.5*exp(-.5*abs(x))};
  \addplot[domain=-4:+4, blue]{0.3*exp(-.5*ln(1 + x^2))};
  \end{axis}
\end{tikzpicture}
\begin{tikzpicture}
  \begin{axis}[
    thick,
    width=.3\textwidth, height=0.9cm,
    xmin=-3,xmax=3,ymin=0,ymax=1.2,
    no markers,
    samples=50,
    axis lines*=left, 
    axis lines*=middle, 
    scale only axis,
    xtick={-1,1},
    xticklabels={},
    ytick={0},
    ] 
\addplot[domain=-3:+3,densely dashed]{.5*x^2};
 \addplot[domain=-3:+3, dashdotted, red]{.5*abs(x)};
  \addplot[domain=-3:+3, blue]{.5*ln(1 + x^2)};
  \end{axis}
\end{tikzpicture}
\begin{tikzpicture}
  \begin{axis}[
    thick,
    width=.3\textwidth, height=0.9cm,
    xmin=-3,xmax=3,ymin=-.8,ymax=.8,
    no markers,
    samples=50,
    axis lines*=left, 
    axis lines*=middle, 
    scale only axis,
    xtick={-1,1},
    xticklabels={},
    ytick={0},
    ] 
\addplot[domain=-3:3,densely dashed]{x};
 \addplot[domain=-3:0, dashdotted, red]{-.5};
  \addplot[domain=0:3, dashdotted, blue]{.5};
\addplot[domain=-3:3]{x/(1 + x^2)};
  \end{axis}
\end{tikzpicture}
    \caption{\label{GLT-KF}
{Gaussian (dashed black line), Laplace (dash-dotted red line), and Student's $t$ (solid blue line); Densities (left plot), Negative Log Likelihoods (center plot), and Influence Functions (right plot).}
   {Student's $t$-density has heavy tails, a non-convex log-likelihood, and a re-descending influence function.}}
    
\end{figure}
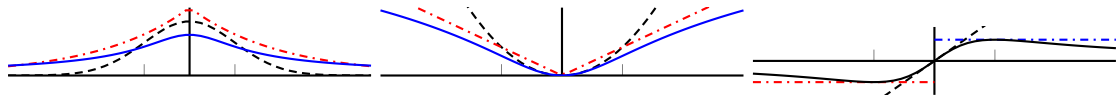

\section{Formulation Examples}

\subsection{Robust Penalties}
In statistics, robust penalties are used since they penalize large residuals less than the classic quadratic norm. The 
Huber penalty is quadratic at the origin and becomes linear for large values of the parameter. 
Student's $t$ based penalties are {\it sublinear} and have been shown to be extremely robust against outlier noise
in the seismic context~\citep{AravkinFHV:2012}. 
Practical implementation simply requires replacement of the least-squares penalty by one of the form $\sum_i\log(\nu + r_i^2)$, 
where the sum runs across a residual vector. The parameter $\nu$ can be tuned by cross validation, or using 
a generalized statistical formulation~\citep{AravkinVanLeeuwen2012}.

\subsection{Sparse Regularization and Constraints}
A great deal of work has underlined the importance of sparse regularization, both theoretically~\citep{Donoho2006_CS,candes2006nos}
and in practice in application to seismic problems (see e.g.,~\cite{Herrmann11TRfcd} and references therein). 
Curvelets~\citep{Demanet2006td} have been a particularly popular choice of transform space, though recent work also advocates for adaptive learning of dictionaries~\citep{horesh2009sensitivity}. Once an efficient transform is found, $\ell_1$-norm solvers offer improved recovery of the parameters of interest. 
The non-smoothness of the regularizer is a key feature and cannot be circumvented.

\section{Results}
\subsection {Spectral Elements Time Domain 2D Model}
\begin{figure}
        \centering

         \begin{subfigure}[t]{0.3\textwidth}
                \includegraphics[width=\textwidth]{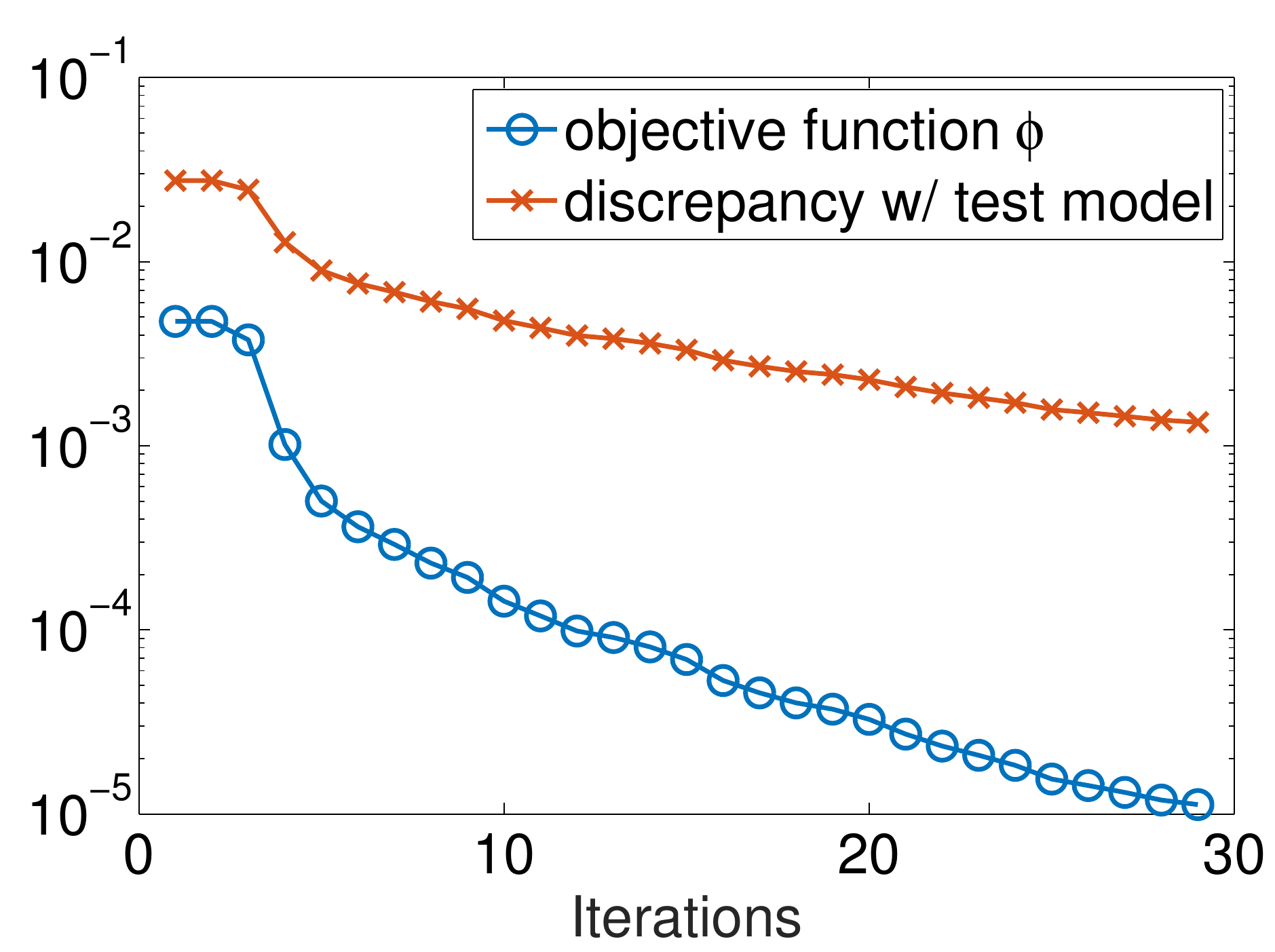}
                \caption{$\rho=\rho_{LS}$, $R=0$}
                \label{fig:A}
        \end{subfigure}
        \begin{subfigure}[t]{0.3\textwidth}
                \includegraphics[width=\textwidth]{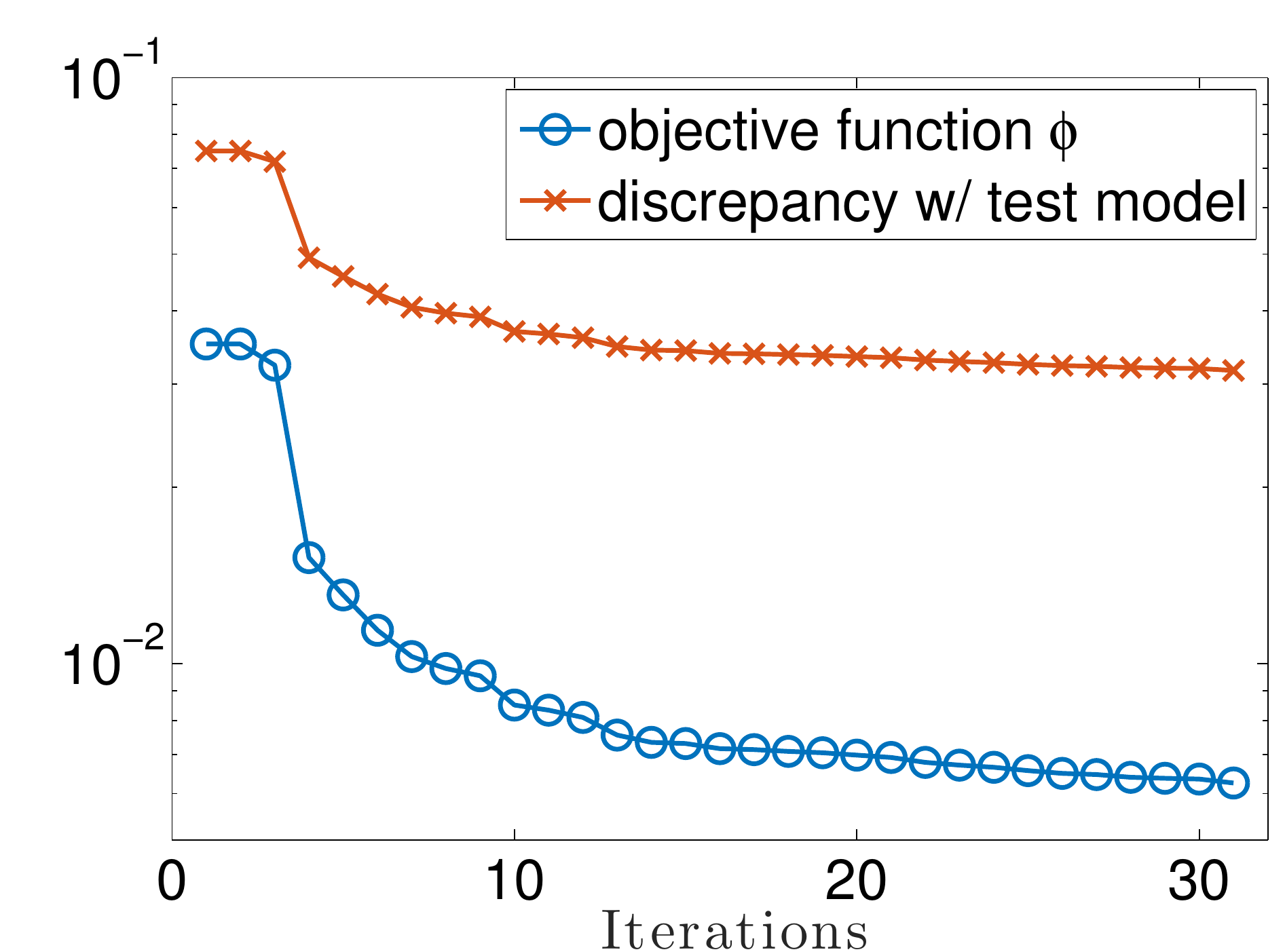}
                \caption{$\rho=\rho_{LS}$, $R=R_{lasso}$}
                \label{fig:B}
        \end{subfigure}
                \begin{subfigure}[t]{0.3\textwidth}
                \includegraphics[width=\textwidth]{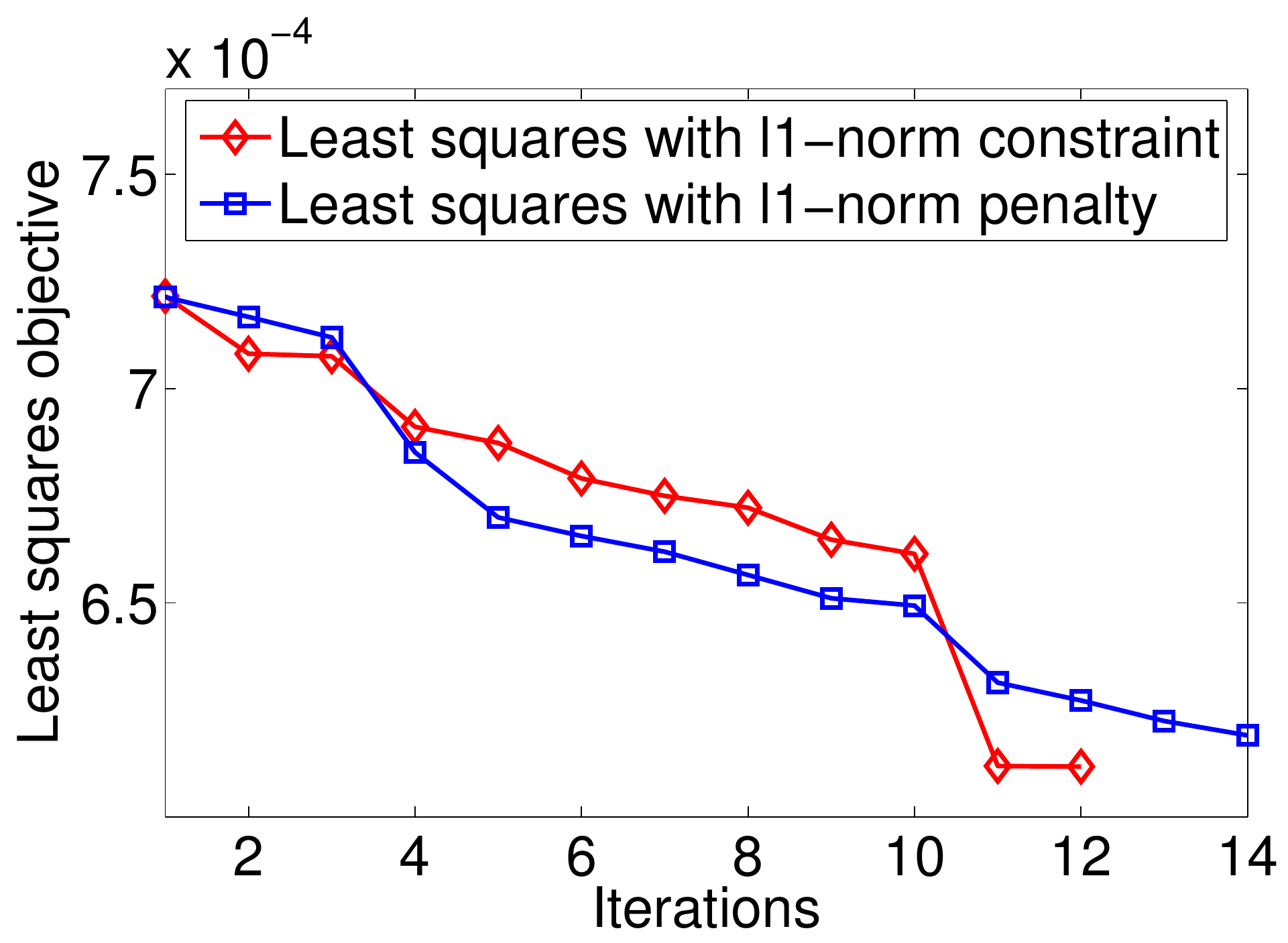}
                \caption{$\rho=\rho_{LS}$, $R=R_{lasso}$ (blue),
                $R=R_{\ell_1-ball}$ (red)}
                \label{fig:E}
        \end{subfigure}

        \begin{subfigure}[t]{0.3\textwidth}
                \includegraphics[width=\textwidth]{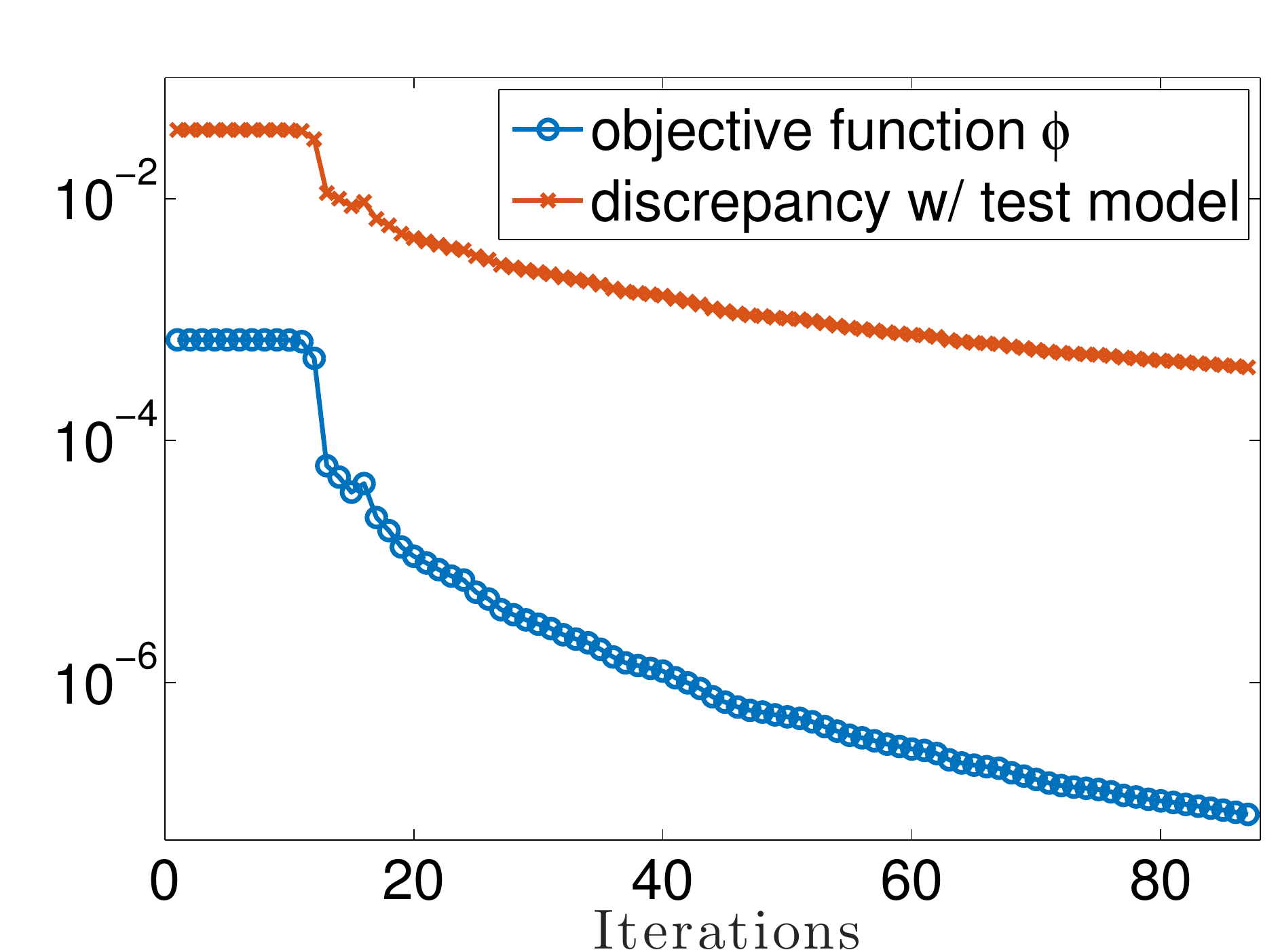}
                \caption{$\rho=\rho_{\text{student}-t}$, $R=0$}
                \label{fig:C}
        \end{subfigure}
        \begin{subfigure}[t]{0.3\textwidth}
                \includegraphics[width=\textwidth]{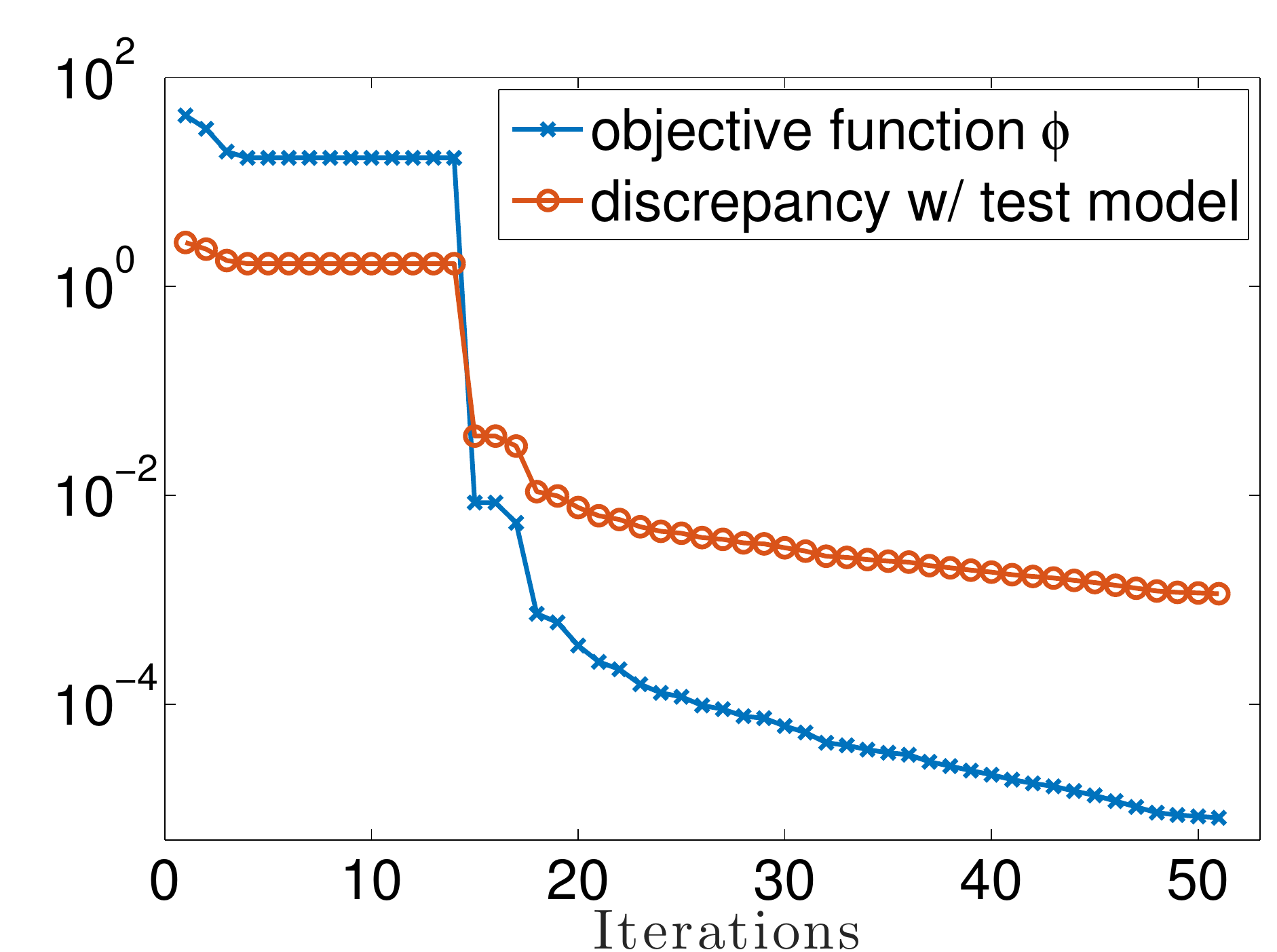}
                \caption{$\rho=\rho_{\text{student}-t}$, $R=R_{lasso}$, $C=$curvelets}
                \label{fig:D}
        \end{subfigure}
        \begin{subfigure}[t]{0.3\textwidth}
                \includegraphics[width=\textwidth]{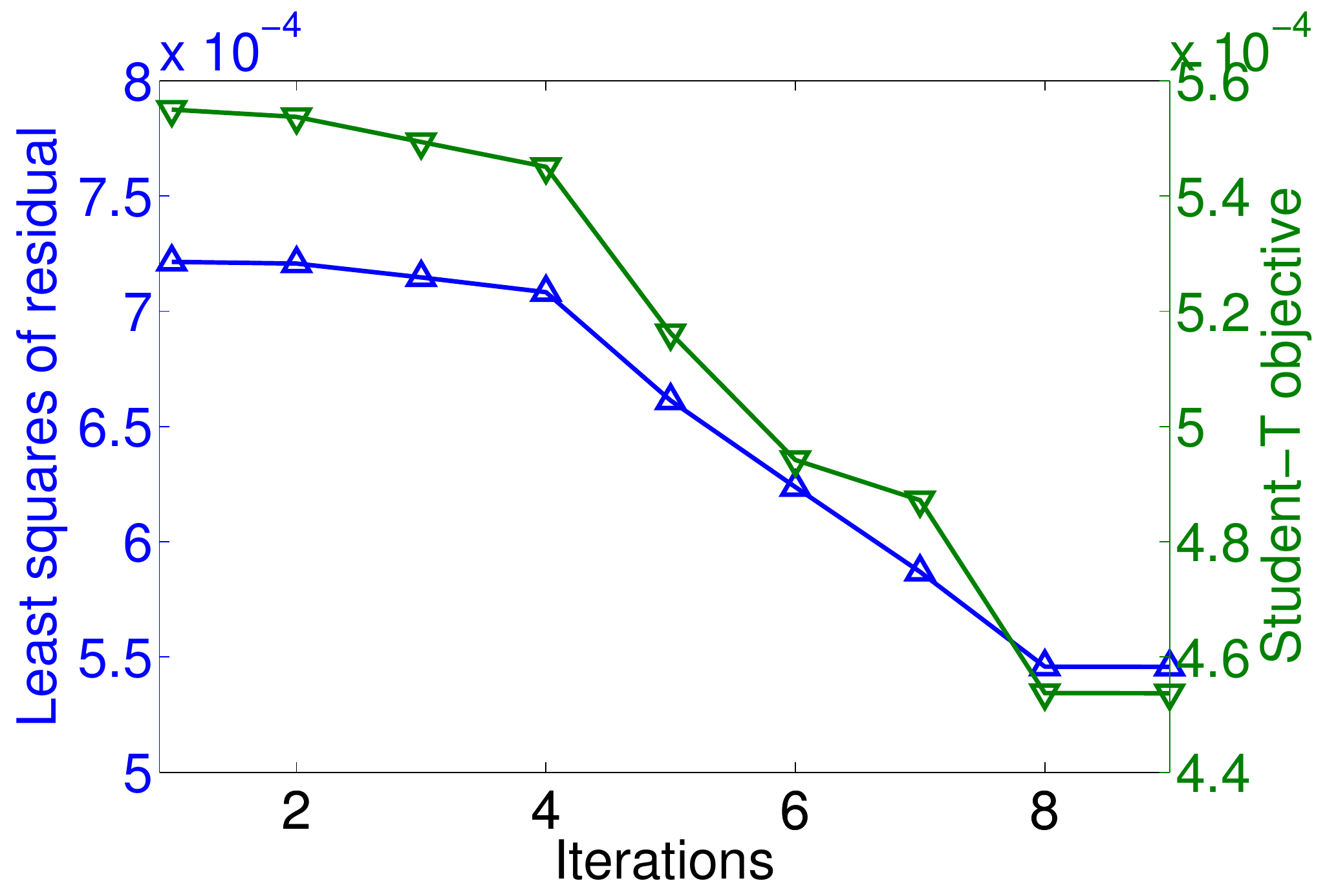}
                \caption{$\rho=\rho_{\text{Student}-t}$, $R=0$}
                \label{fig:F}
        \end{subfigure} 
        
        \caption{Convergence plots, for various $R,\rho,C$. $C=I$ except for (\subref{fig:D})}\label{fig:plots}
\end{figure}

High-order spectral elements in the time domain over a 2D grid is used for simulation, with a test model interpolated from the SEG/EAGE salt dome model using 4 sources and 121 receivers.
The $\ell_2$ ``discrepancy'' is in comparison with the test model, using the Frobenius norm, so for regularized models, it is not necessarily expected to go to zero since the test model may not be 
the optimal solution.
The $1^{st}$ case employs the traditional least-squares loss function (Fig.~\ref{fig:A}). The $2^{nd}$ 
uses least-squares loss with $\ell_1$ penalty (lasso) (Fig.~\ref{fig:B}). 
The $3^{rd}$ test case is Student's $t$ loss function,  $\rho_{student}$, with no regularization functional $R$ (Fig.~\ref{fig:C}). In the $4^{th}$ test case, 
we work in the curvelet domain and use $R_{lasso}$.
 The objective is reduced by 5 and 7 orders of magnitude, resp.\ (for similar initial point).


\subsection{Finite Difference Frequency Domain 3D Model}
Forward simulations are performed in the frequency domain on a 3D grid. An interpolated version of the SEG/EAGE model served as the true model. Data are generated using 6 sources and 60 receivers. 
To validate the code, a standard least-squares inversion is performed, i.e., 
$\rho=\rho_{LS}$, $R=0$, $C=I$.
The objective decreases by 2 orders of magnitude within 6 iterations (results not shown).
The $2^{nd}$ experiment uses the same objective $\rho_{LS}$ (and $C=I$) but adds a constraint
$R$ that bounds the model to be inside the $\ell_1$ ball of radius $\tau$. 
%
$\tau$ can be either set empirically or via cross-validation on a corpus of models.
The initial iterate is 
set to
the standard least-squares solution, so further improvement in the objective value
shown on a linear scale in Fig.~\ref{fig:E} is small.
The $3^{rd}$ experiment is a variant of the second that replaces the $\ell_1$ ball
with an $\ell_1$ penalty of strength $\lambda$, where $\lambda$ is chosen 
empirically but could be again tuned via cross-validation.  The original PQN algorithm 
 does not handle this case since $R$ is
not an indicator function, so we modify the algorithm. In particular, we use
the ``Projected Scaled Sub-gradient + Active Set'' solver~\citep{Schmidt2007a}
as the sub-problem solver.  The convergence of $\rho_{LS}(h(y),D)$ is also shown in
Fig.~\ref{fig:E} in blue
 (but note that we optimized $\rho_{LS}(h(y),D) +
\lambda\|y\|_1$ ).
Finally, we switch from $\rho_{LS}$ to the Student's $t$ based loss
$\rho_{\text{student}-t}$. We show results for $R\equiv 0$ only (Fig.~\ref{fig:F}), but we also have
results for combining with $\ell_1$ ball and $\ell_1$ penalty terms as done
for least-squares.
The left-axis of the plot shows the 
least-squares residual, which decreases, even though this objective was not optimized.


%

\section{Conclusions}
We have presented a general algorithmic framework for large-scale FWI that offers flexible incorporation 
of smooth objectives, as well as various forms of regularization. The framework requires 
only that the proximity operator of the regularization be easily computable, which includes 
projections onto simple sets, as well as sparsity regularization. 

\bibliographystyle{seg}
\bibliography{mybib}

%
%

\end{document}